\documentclass[a4paper,11pt,reqno]{amsart}

\usepackage{dsfont,amssymb}

\chardef\bslash=`\\ 

\newtheorem{theorem}{Theorem}
\newtheorem{corollary}[theorem]{Corollary}
\newtheorem{lemma}[theorem]{Lemma}
\newtheorem{proposition}[theorem]{Proposition}
\theoremstyle{remark}

\theoremstyle{definition}

\newcommand\bp{\begin{proof}}
\newcommand\ep{\end{proof}}

\newcommand\Dhat{{\hat\Delta}}

\newcommand\CC{{\mathcal C}}

\newcommand\E{{\mathcal E}}
\newcommand\F{{\mathcal F}}

\newcommand\RR{{\mathcal R}}
\newcommand\U{{\mathcal U}}

\newcommand\g{{\mathfrak g}}

\newcommand\Aut{\operatorname{Aut}}

\newcommand\End{\operatorname{End}}
\newcommand\Hom{\operatorname{Hom}}

\newcommand\Out{\operatorname{Out}}

\newcommand\tr{\operatorname{tr}}
\newcommand\Vect{\mathcal Vec}

\newcommand{\ad}{\operatorname{ad}}

\newcommand{\C}{{\mathbb C}}

\newcommand\T{{\mathbb T}}
\newcommand\Z{{\mathbb Z}}

\newcommand\enu[1]{\smallskip\newline\makebox[5mm][l]{\rm(#1)}}

\begin{document}

\title{Autoequivalences of the tensor category of $U_q\g$-modules}

\author[S. Neshveyev]{Sergey Neshveyev}
\address{Department of Mathematics, University of Oslo,
P.O. Box 1053 Blindern, NO-0316 Oslo, Norway}

\email{sergeyn@math.uio.no}

\author[L. Tuset]{Lars Tuset}
\address{Faculty of Engineering, Oslo University College,
P.O. Box 4 St.~Olavs plass, NO-0130 Oslo, Norway}
\email{Lars.Tuset@iu.hio.no}

\thanks{Supported by the Research Council of Norway.}

\date{December 21, 2010; revised on July 29, 2011}

\begin{abstract}
We prove that for $q\in\C^*$ not a nontrivial root of unity the cohomology group defined by invariant $2$-cocycles in a completion of $U_q\g$ is isomorphic to $H^2(P/Q;\T)$, where $P$ and $Q$ are the weight and root lattices of $\g$. This implies that the group of autoequivalences of the tensor category of $U_q\g$-modules is the semidirect product of $H^2(P/Q;\T)$ and the automorphism group of the based root datum of $\g$. For $q=1$ we also obtain similar results for all compact connected separable groups.
\end{abstract}

\maketitle

\bigskip

For a tensor category $\CC$  a natural object to study is its group of symmetries, i.e., the group $\Aut^\otimes(\CC)$ of monoidal autoequivalences of $\CC$ identified up to monoidal natural isomorphisms. A more refined version of this group is the tensor category of autoequivalences of $\CC$. It is, for example, used to define what is meant by an action of a group on~$\CC$, which in turn leads to such constructions as equivariantization and crossed products, see e.g.~\cite{Ni} for applications. At the same time there are not many examples for which the group $\Aut^\otimes(\CC)$ is explicitly computed. The aim of this note is to calculate it for the category of representations of the $q$-deformation $G_q$ of a simply connected semisimple compact Lie group $G$. Part of the information about the group of autoequivalences in this case is contained in the work of McMullen~\cite{McMul}, who showed that that the group of automorphisms of the fusion ring of $G$ is isomorphic to $\Out(G)$, that is, to the automorphism group of the based root datum of $\g$. The remaining part is determined by the possible tensor structures one can have on the identity functor, and these are described by the cohomology group defined by invariant $2$-cocycles on the dual $\hat G_q$ of the quantum group $G_q$. Another motivation for computing this cohomology group is the problem of classifying Drinfeld twists that do not necessarily respect braiding; the ones that do respect braiding have been classified in \cite{NT4}.

In a previous paper \cite{NT5} we showed that if $G$ is a compact connected group then the cohomology group defined by invariant unitary $2$-cocycles on $\hat G$ is isomorphic to $H^2(\widehat{Z(G)};\T)$ and we conjectured that for semisimple Lie groups a similar result holds for the $q$-deformation of $G$. We will prove that this is indeed the case using techniques from our earlier paper~\cite{NT4}, where we considered symmetric cocycles and were inspired by the proof of Kazhdan and Lusztig of the equivalence of the Drinfeld category and the category of $U_q\g$-modules~\cite{KL1}. For $q=1$ this gives an alternative proof of the main results in~\cite[Section~2]{NT5} and allows us, at least in the separable case, to extend those results to non-unitary cocycles relying neither on ergodic actions nor on reconstruction theorems. At the same time this proof is less transparent than that in~\cite{NT5} and, as opposed to~\cite{NT5}, relies heavily on the structure and representation theory of compact Lie groups.

\bigskip

We will follow the notation and conventions in~\cite{NT4}. Let $G$ be a simply connected semisimple compact Lie group, $\g$ its complexified Lie algebra,  $q\in\C^*$ not a nontrivial root of unity. Fix a Cartan subalgebra of $\g$ and a system $\{\alpha_1,\dots,\alpha_r\}$ of simple roots. The weight and root lattices are denoted by~$P$ and~$Q$, respectively. For weight $\lambda\in P$ denote by $\lambda(i)$ the coefficients of $\lambda$ in the basis consisting of fundamental weights. Take the $\ad$-invariant symmetric form on $\g$ such that $(\alpha,\alpha)=2$ for every short root in every simple component of $\g$, and put $d_i=(\alpha_i,\alpha_i)/2$ and $q_i=q^{d_i}$.\footnote[1]{Our main result, Theorem~\ref{tmain}, is valid for any $\ad$-invariant symmetric form on $\g$ such that its restriction to the real Lie algebra of $G$ is negative definite, under the assumption that either $q=1$ (in which case the choice of a form does not matter) or that $q_i$ is not a root of unity for all $i$.} For $q\ne1$ consider the quantized universal
enveloping algebra $U_q\g$ with generators $E_i$, $F_i$ and $K_i$, $1\le i\le r$, so that we in particular have $$K_iE_jK_i^{-1}=q_i^{a_{ij}}E_j,\ \ K_iF_jK_i^{-1}=q_i^{-a_{ij}}F_j,\ \ E_iF_j-F_jE_i=\delta_{ij}(K_i-K_i^{-1})/(q_i-q_i^{-1}).$$ Recall that a $U_q\g$-module $V$ is called admissible if $V=\oplus_{\lambda\in P}V(\lambda)$, where $V(\lambda)$ consists of vectors~$v\in V$ such that $K_iv=q_i^{\lambda(i)}v$ for all $i$. Denote by $\CC_q(\g)$ the tensor category of admissible finite dimensional $U_q\g$-modules. For $q=1$ denote by $\CC(\g)=\CC_1(\g)$ the usual tensor category of finite dimensional $U\g$-modules. Let $\U(G_q)$ be the endomorphism ring of the forgetful functor $\CC_q(\g)\to\Vect$. If for every dominant integral weight $\mu\in P_+$ we fix an irreducible $U_q\g$-module $V_\mu$ with highest weight~$\mu$, then the ring $\U(G_q)$ can be identified with $\prod_{\mu\in P_+}\End(V_\mu)$. The comultiplication on $U_q\g$ extends to a homomorphism $\Dhat_q\colon\U(G_q)\to\U(G_q\times G_q)=\prod_{\mu,\eta\in P_+}\End(V_\mu\otimes V_\eta)$.

\smallskip

An invertible element $\E\in \U(G_q\times G_q)$ is called a $2$-cocycle on $\hat G_q$ if
$$
(\E\otimes1)(\Dhat_q\otimes\iota)(\E)=(1\otimes\E)(\iota\otimes\Dhat_q)(\E).
$$
A cocycle is called invariant if it commutes with elements in the image of $\Dhat_q$. The set of invariant $2$-cocycles
forms a group under multiplication, which we denote by $Z^2_{G_q}(\hat G_q;\C^*)$. Cocycles of the form $(a\otimes a)\Dhat_q(a)^{-1}$, where $a$ is an invertible element in the center of $\U(G_q)$, form a subgroup of the center of $Z^2_{G_q}(\hat G_q;\C^*)$. The quotient of $Z^2_{G_q}(\hat G_q;\C^*)$ by this subgroup is denoted by $H^2_{G_q}(\hat G_q;\C^*)$.

\smallskip

The center of $\U(G_q)=\prod_{\mu\in P_+}\End(V_\mu)$ is identified with the algebra of functions on the set $P_+$ of dominant integral weights. By~\cite[Proposition~4.5]{NT4} a function on $P_+$ is a group-like element of~$\U(G_q)$ if and only if it is defined by a character of $P/Q$. Therefore the Hopf algebra of functions on $P/Q$ embeds into the center of~$\U(G_q)$. Hence every $2$-cocycle $c$ on $P/Q$ can be considered as an invariant $2$-cocycle~$\E_c$ on~$\hat G_q$. Explicitly, $\E_c$ acts on $V_\mu\otimes V_\eta$ as multiplication by $c(\mu,\eta)$. We can now formulate our main result.

\begin{theorem} \label{tmain}
The homomorphism $c\mapsto \E_c$ induces an isomorphism
$$
H^2(P/Q;\T)\cong H^2_{G_q}(\hat G_q;\C^*).
$$
In particular, if $\g$ is simple and $\g\not\cong\mathfrak{so}_{4n}(\C)$ then $H^2_{G_q}(\hat G_q;\C^*)$ is trivial, and if $\g\cong\mathfrak{so}_{4n}(\C)$ then $H^2_{G_q}(\hat G_q;\C^*)\cong\Z/2\Z$.
\end{theorem}

The last statement follows from the fact that for simple Lie algebras the group $P/Q$ is cyclic unless $\g\cong\mathfrak{so}_{4n}(\C)$, in which case $P/Q\cong\Z/2\Z\times\Z/2\Z$, see e.g.~Table~IV on page~516 in~\cite{Hel}.

\smallskip

Note that for $q>0$ the same result holds for unitary cocycles. This easily follows by polar decomposition,
see~\cite[Lemma~1.1]{NT4}.

\smallskip

In the proof of the theorem we will assume that $q\ne1$, the case $q=1$ is similar and for unitary cocycles is also proved by a different method in~\cite{NT5}.

\smallskip

Our first goal will be to construct a homomorphism $H^2_{G_q}(\hat G_q;\C^*)\to H^2(P/Q;\T)$. For every $\mu\in P_+$ fix a highest weight vector $\xi_\mu\in V_\mu$. Recall~\cite[Section~2]{NT4} that for $\mu,\eta\in P_+$ there exists a unique morphism
$$
T_{\mu,\eta}\colon V_{\mu+\eta}\to V_\mu\otimes V_\eta \ \ \hbox{such that}\ \
\xi_{\mu+\eta}\mapsto\xi_\mu\otimes \xi_\eta.
$$
The image of $T_{\mu,\eta}$ is the isotypic component of $V_\mu\otimes V_\eta$ with highest weight $\mu+\eta$. Hence if~$\E$~is an invariant $2$-cocycle then it acts on this image as multiplication by a nonzero scalar $c_\E(\mu,\eta)$. As in the proof of~\cite[Lemma~2.2]{NT4}, the identity $(T_{\mu,\eta}\otimes\iota)T_{\mu+\eta,\nu}
=(\iota\otimes T_{\eta,\nu})T_{\mu,\eta+\nu}$ immediately implies that $c_\E$ is a $2$-cocycle on $P_+$. Furthermore, the cohomology class $[c_\E]$ of $c_\E$ in $H^2(P_+;\C^*)$ depends only on the class of $\E$ in $H^2_{G_q}(\hat G_q;\C^*)$, since if $a\in \U(G_q)$ is a central element acting on $V_\mu$ as multiplication by a scalar~$a(\mu)$ then the action of $(a\otimes a)\Dhat_q(a)^{-1}$ on the image of $T_{\mu,\eta}$ is multiplication by $a(\mu)a(\eta)a(\mu+\eta)^{-1}$. Thus the map $\E\mapsto c_\E$ defines a homomorphism $H^2_{G_q}(\hat G_q;\C^*)\to H^2(P_+;\C^*)$.

Given a cocycle on $P/Q$, we can consider it as a cocycle on $P$ and then get a cocycle on $P_+$ by restriction. Thus we have a homomorphism $H^2(P/Q;\T)\to H^2(P_+;\C^*)$. It is injective since the quotient map $P_+\to P/Q$ is surjective and a cocycle on $P/Q$ is a coboundary if it is symmetric.

\begin{lemma}
For every invariant $2$-cocycle $\E$ on $\hat G_q$ the class of $c_\E$ in $H^2(P_+;\C^*)$ is contained in the image of $H^2(P/Q;\T)$.
\end{lemma}

\bp Consider the skew-symmetric bi-quasicharacter $b\colon P_+\times P_+\to\C^*$ defined by $$b(\mu,\eta)=c_\E(\mu,\eta)c_\E(\eta,\mu)^{-1}.$$ It extends uniquely to a skew-symmetric bi-quasicharacter on $P$. To prove the lemma it suffices to show that the root lattice $Q$ is contained in the kernel of this extension. Indeed, since $H^2(P/Q;\T)$ is isomorphic to the group of skew-symmetric bi-characters on $P/Q$, it then follows that there exists a cocycle $c$ on $P/Q$ such that the cocycle $c_\E c^{-1}$ on $P_+$ is symmetric. Then by~\cite[Lemma~4.2]{NT3} the cocycle $c_\E c^{-1}$ is a coboundary, so $c_\E$ and the restriction of $c$ to $P_+$ are cohomologous.

To prove that $Q$ is contained in the kernel of $b$, recall~\cite[Section~2]{NT4} that for every simple root $\alpha_i$ and weights $\mu,\eta\in P_+$ with $\mu(i),\eta(i)\ge1$ we can define a morphism
\begin{equation*}
\tau_{i;\mu,\eta}\colon V_{\mu+\eta-\alpha_i}\to V_\mu\otimes
V_\eta\ \ \hbox{such that}\ \ \xi_{\mu+\eta-\alpha_i}\mapsto
[\mu(i)]_{q_i}\xi_\mu\otimes F_i\xi_\eta
-q_i^{\mu(i)}[\eta(i)]_{q_i}F_i\xi_\mu\otimes\xi_\eta.
\end{equation*}
The image of $\tau_{i;\mu,\eta}$ is the isotypic component of $V_\mu\otimes
V_\eta$ with highest weight $\mu+\eta-\alpha_i$. Since the element $\E$ is invariant, it acts on this image as multiplication by a nonzero scalar $c_i(\mu,\eta)$. As in the proof of~\cite[Lemma~2.3]{NT4}, consider now another weight $\nu$ with $\nu(i)\ge1$. The isotypic component of $V_\mu\otimes
V_\eta\otimes V_\nu$ with highest weight $\mu+\eta+\nu-\alpha_i$ has
multiplicity two, and is spanned by the images of $(\iota\otimes
T_{\eta,\nu})\tau_{i;\mu,\eta+\nu}$ and
$(\iota\otimes\tau_{i;\eta,\nu})T_{\mu,\eta+\nu-\alpha_i}$, as well as
by the images of $(T_{\mu,\eta}\otimes\iota)\tau_{i;\mu+\eta,\nu}$ and
$(\tau_{i;\mu,\eta}\otimes\iota)T_{\mu+\eta-\alpha_i,\nu}$. We have
\begin{equation}\label{etauT}
[\eta(i)]_{q_i}(T_{\mu,\eta}\otimes\iota)\tau_{i;\mu+\eta,\nu}
-[\nu(i)]_{q_i}(\tau_{i;\mu,\eta}\otimes\iota)
T_{\mu+\eta-\alpha_i,\nu}=[\mu(i)+\eta(i)]_{q_i}
(\iota\otimes\tau_{i;\eta,\nu})
T_{\mu,\eta+\nu-\alpha_i}.
\end{equation}
Apply the element
$$
\Omega:=(\E\otimes1)(\Dhat_q\otimes\iota)(\E)=(1\otimes\E)(\iota\otimes\Dhat_q)(\E)
$$
to this identity. The morphisms $(T_{\mu,\eta}\otimes\iota)\tau_{i;\mu+\eta,\nu}$, $(\tau_{i;\mu,\eta}\otimes\iota)
T_{\mu+\eta-\alpha_i,\nu}$ and $(\iota\otimes\tau_{i;\eta,\nu})
T_{\mu,\eta+\nu-\alpha_i}$ are eigenvectors of the operator of multiplication by $\Omega$ on the left with eigenvalues $c_\E(\mu,\eta)c_i(\mu+\eta,\nu)$, $c_i(\mu,\eta)c_\E(\mu+\eta-\alpha_i,\nu)$ and $c_i(\eta,\nu)c_\E(\mu,\eta+\nu-\alpha_i)$, respectively. Since the morphisms  $(T_{\mu,\eta}\otimes\iota)\tau_{i;\mu+\eta,\nu}$ and
$(\tau_{i;\mu,\eta}\otimes\iota)T_{\mu+\eta-\alpha_i,\nu}$ are linearly independent, by applying $\Omega$ to \eqref{etauT} we conclude that these three eigenvalues coincide. In particular,
$$c_i(\mu,\eta)c_\E(\mu+\eta-\alpha_i,\nu)=c_i(\eta,\nu)c_\E(\mu,\eta+\nu-\alpha_i).$$
Applying this to $\eta=\nu=\mu$ we get
$$
b(2\mu-\alpha_i,\mu)=1.
$$
Since $b$ is skew-symmetric, this gives $b(\alpha_i,\mu)=1$. The latter identity holds for all $\mu\in P_+$ with $\mu(i)\ge1$. Since every element in $P$ can be written as a difference of two such elements $\mu$, it follows that $\alpha_i$ is contained in the kernel of $b$.
\ep

Therefore the map $\E\mapsto c_\E$ induces a homomorphism $H^2_{G_q}(\hat G_q;\C^*)\to H^2(P/Q;\T)$. Clearly, it is a left inverse of the homomorphism $H^2(P/Q;\T)\to H^2_{G_q}(\hat G_q;\C^*)$, $[c]\mapsto[\E_c]$, constructed earlier. Thus it remains to prove that the homomorphism $H^2_{G_q}(\hat G_q;\C^*)\to H^2(P/Q;\T)$ is injective.

\smallskip

Assume that $\E$ is an invariant $2$-cocycle such that the cocycle $c_\E$ on $P_+$ is a coboundary. Our goal is to show that $\E$ is the coboundary of a central element in $\U(G_q)$. We will follow the strategy in \cite{NT4}, where this was shown under the additional assumption that $\E$ is symmetric, that is, $\RR_\hbar\E=\E_{21}\RR_\hbar$ for an $R$-matrix $\RR_\hbar\in\U(G_q\times G_q)$, which depends on the choice of a number $\hbar\in\C$ such that $q=e^{\pi i \hbar}$.

The first step in \cite{NT4}, see the discussion following Lemma~2.2 in \cite{NT4}, was to show that $\E$ is cohomologous to a cocycle such that
\begin{equation} \label{econd}
\E T_{\mu,\eta}=T_{\mu,\eta}\ \ \hbox{and}\ \ \E \tau_{i;\nu,\omega}=\tau_{i;\nu,\omega}
\end{equation}
for all $\mu,\eta\in P_+$, $1\le i\le r$ and $\nu,\omega\in P_+$ such that $\nu(i),\omega(i)\ge1$. This part goes through in the non-symmetic case without any changes, as the symmetry requirement was needed only to conclude that $c_\E$ is a coboundary.

Therefore to prove the injectivity of $H^2_{G_q}(\hat G_q;\C^*)\to H^2(P/Q;\T)$  it suffices to establish the following result, which extends \cite[Corollary~4.4]{NT4}.

\begin{proposition} \label{pcond}
If $\E$ is an invariant $2$-cocycle on $\hat G_q$  with property \eqref{econd} then $\E=1$.
\end{proposition}

The proof of this statement in \cite{NT4} for symmetric cocycles is based on considering the action of $\E$ on a comonoid representing the forgetful functor on $\CC_q(\g)$. Recall briefly how this comonoid, essentially constructed by Kazhdan and Lusztig, is defined. For every weight $\mu\in P_+$ fix an irreducible $U_q\g$-module $\bar V_\mu$ with lowest weight $-\mu$ and a lowest weight vector $\bar\xi_\mu$. For $\lambda\in P$ and $\mu,\eta\in P_+$ such that $\lambda+\mu\in P_+$, there exists a unique morphism
$$
\tr^{\eta}_{\mu,\lambda+\mu}\colon \bar V_{\mu+\eta}\otimes
V_{\lambda+\mu+\eta}\to \bar V_{\mu}\otimes V_{\lambda+\mu}\ \ \hbox{such that}\ \
\bar\xi_{\mu+\eta}\otimes\xi_{\lambda+\mu+\eta}
\mapsto\bar\xi_{\mu}\otimes\xi_{\lambda+\mu}.
$$
Using these morphisms define an inverse limit $U_q\g$-module
$$
M_\lambda=\lim_{\xleftarrow[\mu]{}}\bar V_\mu\otimes
V_{\lambda+\mu}.
$$
Denote by $\tr_{\mu,\lambda+\mu}$ the canonical map $M_\lambda\to \bar V_\mu\otimes
V_{\lambda+\mu}$. The module $M_\lambda$ is considered as a topological $U_q\g$-module with a base of
neighborhoods of zero formed by the kernels of the maps $\tr_{\mu,\lambda+\mu}$, while all modules in our category $\CC_q(\g)$ are considered with discrete topology. Then $\Hom_{U_q\g}(M_\lambda,V)$ is the inductive limit of the spaces
$\Hom_{U_q\g}(\bar V_{\mu}\otimes V_{\lambda+\mu},V)$. The vectors $\bar\xi_\mu\otimes\xi_{\lambda+\mu}$ define a topologically cyclic vector $\Omega_\lambda\in M_\lambda$. For any finite dimensional admissible $U_q\g$-module $V$ the map
$$
\eta_V\colon\Hom_{U_q\g}(\oplus_\lambda M_\lambda,V)\to V,\ \ \eta_V(f)=\sum_\lambda f(\Omega_\lambda),
$$
is an isomorphism, so the topological $U_q\g$-module $M=\oplus_\lambda M_\lambda$ represents the forgetful functor. Furthermore, the representation of $U_q\g$ in the endomorphism ring of the forgetful functor is implemented by the antihomomorphism $\pi\colon U_q\g\to\End_{U_q\g}(M)$ defined by $\pi(E_i)\Omega_\lambda=E_i\Omega_{\lambda-\alpha_i}$, $\pi(F_i)\Omega_i=F_i\Omega_{\lambda+\alpha_i}$ and $\pi(K_i)\Omega_{\lambda}=q_i^{\lambda(i)}\Omega_\lambda$. In other words, $M$ is a $U_q\g$-bimodule.

It was shown in \cite[Section~4]{NT4}, see the arguments up to (but not including) Lemma~4.3 there, that  condition \eqref{econd} is exactly what is needed to define an action of any invariant cocycle $\E$ satisfying~\eqref{econd} on the $U_q\g$-bimodule $M$. More precisely, we showed that there exist a character $\chi$ of $P/Q$, an invertible morphism $\E_0$ of $M=\oplus_\lambda M_\lambda$ onto itself preserving the direct sum decomposition, and an invertible element~$c$ in the center of $\U(G_q)$ such that
\begin{equation} \label{eresult}
\tr_{\mu,\lambda+\mu}\E_0=\chi(\mu)^{-1}\E\tr_{\mu,\lambda+\mu}\ \ \hbox{and}\ \ \eta_V(f\E_0)=c\,\eta_V(f)
\end{equation}
for all $\mu\in P_+$, $\lambda\in P$ such that $\lambda+\mu\in P_+$, and for all finite dimensional admissible $U_q\g$-modules~$V$ and $f\in\Hom_{U_q\g}(M_\lambda,V)$. We will show now that this is already enough to conclude that $\E$ is, in fact, symmetric.

\bp[Proof of Proposition~\ref{pcond}] We want to show that $\RR_\hbar\E=\E_{21}\RR_\hbar$ for some $\hbar$ such that $q=e^{\pi i \hbar}$. We will prove a stronger statement: $\sigma \E=\E\sigma$ for any braiding $\sigma$ on $\CC_q(\g)$.

By \eqref{eresult}, since $\tr_{\mu,\lambda+\mu}(\Omega_\lambda)=\bar\xi_\mu\otimes\xi_{\lambda+\mu}$, for any $\mu,\eta,\nu\in P_+$ and $f\in\Hom_{U_q\g}(\bar V_\mu\otimes V_\eta,V_\nu)$ we have
$$
\chi(\mu)^{-1}f\E(\bar\xi_\mu\otimes\xi_\eta)=c(\nu)f(\bar\xi_\mu\otimes\xi_\eta).
$$
As the vector $\bar\xi_\mu\otimes\xi_\eta$ is cyclic, this means that $f\E=\chi(\mu)c(\nu)f$. Since this is true for all $f$, we conclude that $\E$ acts on the isotypic component of $\bar V_\mu\otimes V_\eta$ with highest weight $\nu$ as multiplication by $\chi(\mu)c(\nu)$. In other words, $\E$ acts on the isotypic component of $V_\mu\otimes V_\eta$ with highest weight $\nu$ as multiplication by $\chi(\bar\mu)c(\nu)$. It follows that
$$
\sigma\E=\chi(\bar\mu-\bar\eta)\E\sigma\ \ \hbox{on}\ \ V_\mu\otimes V_\eta.
$$
But by assumption \eqref{econd} the element $\E$ is the identity on the isotypic component of $V_\mu\otimes V_\eta$ with highest weight $\mu+\eta$, so by considering the above identity on this isotypic component we conclude that $\chi(\bar\mu-\bar\eta)=1$. Thus $\chi$ is the trivial character and  $\sigma\E=\E\sigma$. By \cite[Corollary~4.4]{NT4} we then get that $\E=1$. This completes the proof of Proposition~\ref{pcond} and hence of Theorem~\ref{tmain}.
\ep

As our first application we will classify Drinfeld twists, relating the coproducts on $U_q\g$ and $U\g$, that do not necessarily respect braiding.

\begin{corollary}
Let $\varphi\colon \U(G_q)\to\U(G)$ be an isomorphism extending the canonical identifications of the centers of these algebras with the algebra of functions on $P_+$, and let $\hbar$ be such that $q=e^{\pi i\hbar}$. Suppose $\F\in \U(G\times G)$ is an invertible element such that
\enu{i} $(\varphi\otimes\varphi)\Dhat_q=\F\Dhat\varphi(\cdot)\F^{-1}$;
\enu{ii} the element $(\iota\otimes\Dhat)(\F^{-1}) (1\otimes\F^{-1})(\F\otimes1)(\Dhat\otimes\iota)(\F)$ coincides with Drinfeld's KZ-associator $\Phi_{KZ}(\hbar t_{12},\hbar t_{23})$, where $t\in\g\otimes\g$ is the element defined by our fixed $\ad$-invariant form on $\g$.

\smallskip

Assume $\F'\in \U(G\times G)$ is another element with the same properties. Then there exist a $\T$-valued $2$-cocycle $c$ on $P/Q$ and an invertible central element $a\in \U(G)$ such that $\F'=\E_c\F(a\otimes a)\Dhat(a)^{-1}$.
\end{corollary}

\bp The proof is similar to that of \cite[Theorem~5.2]{NT4}. Define $\E=(\varphi^{-1}\otimes\varphi^{-1})(\F'\F^{-1})\in\U(G_q\times G_q)$. It is easy to check that $\E$ is an invariant $2$-cocycle on $\hat G_q$. By Theorem~\ref{tmain}, $\E=\E_c(b\otimes b)\Dhat_q(b)^{-1}$ for a $2$-cocycle $c$ on $P/Q$ and a central element $b\in\U(G_q)$. Letting $a=\varphi(b)$, we obtain
$
\F'=\E_c(a\otimes a)(\varphi\otimes\varphi)(\Dhat_q(b)^{-1})\F=\E_c\F(a\otimes a)\Dhat(a)^{-1}$.
\ep

Note that this corollary implies that the Dirac operator defined as in~\cite{NT2} is the same (for fixed~$\varphi$) for any choice of a unitary element $\F$ satisfying properties (i) and (ii). This extends \cite[Theorem~6.1]{NT4}.

\smallskip

We now turn to our main application, the computation of the group of $\C$-linear monoidal autoequivalences of $\CC_q(\g)$ identified up to monoidal natural isomorphisms.

Any automorphism $\alpha$ of the based root datum $\Psi_\g$ of $\g$ defines an automorphism of the Hopf algebra $U_q\g$, hence an autoequivalence~$\tilde\alpha$ of~$\CC_q(\g)$. On the other hand, for any $2$-cocycle $c$ on $P/Q$ we can define an autoequivalence~$\beta_c$ which acts trivially on objects and morphisms, while the tensor structure is given by the action of~$\E_c^{-1}$. It turns out that any autoequivalence of $\CC_q(\g)$ is monodially naturally isomorphic to a composition of two autoequivalences defined either by an automorphism of $\Psi_\g$ or by a cocycle on $P/Q$.

\begin{theorem} \label{tequiv}
The group of $\C$-linear monoidal autoequivalences of the tensor category $\CC_q(\g)$ is canonically isomorphic to $H^2(P/Q;\T)\rtimes\Aut(\Psi_\g)$.
\end{theorem}

\bp The proof is essentially identical to~\cite[Theorem~2.5]{NT5}. Briefly, by a result of McMullen~\cite{McMul} any automorphism of the fusion ring of $\CC_q(\g)$, mapping irreducibles into irreducibles, is implemented by an automorphism of $\Psi_\g$. Hence for any autoequivalence $\gamma$ of $\CC_q(\g)$ there exists a unique automorphism~$\alpha$ of~$\Psi_\g$ such that $\tilde\alpha\circ\gamma$ maps every object into an isomorphic one; that is, ignoring the tensor structure, $\tilde\alpha\circ\gamma$ is naturally isomorphic to the identity functor. Possible tensor structures on the identity functor are, in turn, described by invariant $2$-cocycles on $\hat G_q$.
\ep

We next consider $q=1$ and extend the above results to compact connected groups.

The group $P/Q$ is canonically identified with the dual of the center $Z(G)$ of the group $G$, and so, for $q=1$, Theorem~\ref{tmain} can be formulated as $H^2_G(\hat G;\C^*)\cong H^2(\widehat{Z(G)};\C^*)$.

\begin{theorem}
For any compact connected separable group $G$ we have a canonical isomorphism
$$
H^2_G(\hat G;\C^*)\cong H^2(\widehat{Z(G)};\C^*).
$$
\end{theorem}

\bp For Lie groups the proof is essentially the same as above, with $P$ replaced by the weight lattice of a maximal torus of $G$. In the general case we have a homomorphism $H^2(\widehat{Z(G)};\C^*)\to H^2_G(\hat G;\C^*)$ obtained by considering $\U(Z(G))$ as a subring of $\U(G)$. To construct the inverse homomorphism, for every quotient $H$ of $G$ which is a Lie group consider the composition
$$
H^2_G(\hat G;\C^*)\to H^2_H(\hat H;\C^*)\to H^2(\widehat{Z(H)};\C^*),
$$
where the first homomorphism is defined using the quotient map $\U(G)\to\U(H)$. The map $Z(G)\to Z(H)$ is surjective (since this is true for Lie groups), so $Z(G)$ is the inverse limit of the groups $Z(H)$. Then $H^2(\widehat{Z(G)};\C^*)$ is the inverse limit of the groups $H^2(\widehat{Z(H)};\C^*)$. Therefore the above maps
$H^2_G(\hat G;\C^*)\to H^2(\widehat{Z(H)};\C^*)$ define a homomorphism $H^2_G(\hat G;\C^*)\to H^2(\widehat{Z(G)};\C^*)$. It is clearly a left inverse of the map $H^2(\widehat{Z(G)};\C^*)\to H^2_G(\hat G;\C^*)$, so it remains to show that it is injective.

In other words, we have to check that if $\E$ is an invariant cocycle on $\hat G$ such that its image in~$\U(H\times H)$ is a coboundary for every Lie group quotient $H$ of $G$, then $\E$ itself is a coboundary.
If~$\E$ were unitary, this could be easily shown by taking a weak operator limit point of cochains, see the proof of~\cite[Theorem~2.2]{NT5}, and would not require the separability of $G$. In the non-unitary case we can argue as follows.

Since $G$ is separable, there exists a decreasing sequence of closed normal subgroups $N_n$ of $G$ such that $\cap_{n\ge1} N_n=\{e\}$ and the quotients $H_n=G/N_n$ are Lie groups. Let $\E_n$ be the image of~$\E$ in $\U(H_n\times H_n)$. By assumption there exist invertible central elements $c_n\in \U(H_n)$ such that $\E_n=(c_n\otimes c_n)\Dhat(c_n)^{-1}$. For a fixed $n$ consider the image $a$ of $c_{n+1}$ in $\U(H_n)$. Then $c_na^{-1}$ is a central group-like element in $\U(H_n)$. By~\cite[Theorem~A.1]{NT4} it is therefore defined by an element of the center of the complexification $(H_n)_\C$ of $H_n$. Since the homomorphism $(H_{n+1})_\C\to (H_n)_\C$ is surjective, we conclude that there exists a central group-like element $b$ in $\U(H_{n+1})$ such that its image in $\U(H_n)$ is $c_na^{-1}$. Replacing $c_{n+1}$ by $c_{n+1}b$ we get an element such that $\E_{n+1}=(c_{n+1}\otimes c_{n+1})\Dhat(c_{n+1})^{-1}$  and the image of $c_{n+1}$ in $\U(H_n)$ is $c_n$. Applying this procedure inductively we can therefore assume that the image of $c_{n+1}$ in $\U(H_n)$ is $c_n$ for all $n\ge1$. Then the elements $c_n$ define a central element $c\in \U(G)$ such that $\E=(c\otimes c)\Dhat(c)^{-1}$.
\ep

In~\cite[Theorem~2.5]{NT5} we computed the group of autoequivalences of the C$^*$-tensor category of finite dimensional unitary representations of $G$. The above theorem and the same arguments as in the proof of Theorem~\ref{tequiv} allow us to get a similar result ignoring the C$^*$-structure.

\begin{theorem}
For any compact connected separable group $G$, the group of $\C$-linear monoidal autoequivalences of the category of finite dimensional representations of $G$ is canonically isomorphic to $H^2(\widehat{Z(G)};\C^*)\rtimes\Out(G)$.
\end{theorem}



\bigskip

\begin{thebibliography}{99}

\bibitem{Hel}
S. Helgason, {\em Differential geometry, Lie groups, and symmetric spaces}, Graduate Studies in Mathematics, {\bf 34}, American Mathematical Society, Providence, RI, 2001.

\bibitem{KL1}
D. Kazhdan and G. Lusztig, {\em Tensor structures arising from affine
Lie algebras. III}, J. Amer. Math. Soc. {\bf 7} (1994), 335--381.

\bibitem{McMul}
J.R. McMullen, {\em
On the dual object of a compact connected group}, Math. Z. {\bf 185} (1984), 539--552.

\bibitem{NT2}
S. Neshveyev and L. Tuset, {\em The Dirac operator on compact quantum
groups},  J. Reine Angew. Math. {\bf 641} (2010), 1--20.

\bibitem{NT4}
S. Neshveyev and L. Tuset, {\em Symmetric invariant cocycles on the duals of $q$-deformations}, Adv. Math. {\bf 227} (2011), 146--169.

\bibitem{NT3} S. Neshveyev and L. Tuset, {\em Notes on the Kazhdan-Lusztig
    theorem on equivalence of the Drinfeld category and the category of
    $U_q(\g)$-modules}, preprint arXiv: 0711.4302v1 [math.QA].

\bibitem{NT5}
S. Neshveyev and L. Tuset, {\em On second cohomology of duals of compact groups}, preprint arXiv: 1011.4569v5 [math.OA].

\bibitem{Ni}
D. Nikshych, {\em Non-group-theoretical semisimple Hopf algebras from group actions on fusion categories}, Selecta Math. (N.S.) {\bf 14}  (2008), 145--161.

\end{thebibliography}
\end{document}